\magnification\magstep1
\baselineskip = 18pt
\def\n{\noindent}
\def\qed{{\hfill{\vrule height7pt width7pt
depth0pt}\par\bigskip}}
\overfullrule = 0pt

\def\pf{\n {\bf Proof.\ }}
\centerline{\bf Complex Interpolation}
\centerline{\bf and} 
\centerline{\bf  Regular Operators Between Banach
Lattices} \centerline{by}
\centerline{{ Gilles Pisier}\footnote*{Supported in
part by N.S.F. grant DMS 9003550}}
\vskip24pt
 \centerline{Texas A. and M. University}

\centerline{College Station, TX 77843, U. S. A.}

\centerline{and}

\centerline{Universit\'e Paris 6}

\centerline{Equipe d'Analyse, Bo\^\i te 186,
75252 Paris Cedex 05, France}

 \bigskip	\bigskip

{\bf Abstract.} {  We study certain interpolation and
extension properties of the space of regular operators
between two Banach lattices. Let $R_p$ be the space of all
the regular (or equivalently order bounded) operators on
$L_p$ equipped with the regular norm. We prove the
isometric identity $R_p = (R_\infty,R_1)^\theta$ if
$\theta = 1/p$, which shows that the spaces $(R_p)$ form
an interpolation scale relative to Calder\'on's
interpolation method. We also prove that if $S\subset L_p$
is a subspace, every regular operator $u : S \to L_p$
admits a regular extension $\tilde u : L_p \to L_p$ with
the same regular norm. This extends a result due to
Mireille L\'evy in the case $p = 1$.  Finally, we apply
these ideas to the Hardy space $H^p$ viewed as a subspace
of $L_p$ on the circle. We show that the space of regular
operators from $H^p$ to $L_p$ possesses a similar
interpolation property as the spaces $R_p$ defined above.}

\vfill\eject
In a recent paper [HeP] we have observed that the real interpolation spaces
associated to the couples
$$(B(c_0,\ell_\infty), \quad B(\ell_1,\ell_1))$$
and
$$(B(L_\infty,L_\infty), \quad B(L_1,L_1))$$
can be described and an equivalent of the $K_t$-functional
can be given (cf. [HeP]). It is natural to wonder whether
analogous results hold for the complex interpolation
method and this is the subject of the present paper.

Let $X_0,X_1$ be Banach lattices of measurable functions defined on some
set $S$ (we are deliberately vague, see e.g. [LT] for a detailed theory).
We will denote by
$$X_\theta = X^{1-\theta}_0 X_1^\theta$$
the space of all measurable functions $f$ on the set $S$ such that there
are $f_0 \in X_0$, $f_1\in X_1$ satisfying $|f| \le
|f_0|^{1-\theta} |f_1|^\theta$ and we let
$$\|f\|_{X_\theta} = \inf\{\|f_0\|^{1-\theta}_{X_0}
\|f_1\|^\theta_{X_1}\}$$
where the infimum runs over all possible   decompositions
of $f$.

We will denote by $(X_0,X_1)_\theta$ and $(X_0,X_1)^\theta$ the complex
interpolation spaces as defined for examples in [BL]. Recall the
fundamental identity (due to Calder\'on)
$$X^{1-\theta}_0 X_1^{\theta} = (X_0,X_1)^\theta$$
with identical norms, which is valid under the assumption
that the unit ball of $X^{1-\theta}_0 X_1^{\theta}$ is
closed in $X_0+X_1$ (see [Ca]). Moreover, if either $X_0$
or $X_1$ is reflexive (cf. also [HP] and [B]) we have
$$(X_0,X_1)_\theta = (X_0,X_1)^\theta,$$ with identical
norms. In particular when $X_0,X_1$ are finite dimensional
spaces, there is no need to distinguish $(X_0,X_1)_\theta$
and $(X_0,X_1)^\theta$. We will use this fact repeatedly in
the sequel. Let $c_0$ (resp. $\ell_\infty$) be the space
of all sequences of complex scalars tending to zero (resp.
bounded) at infinity equipped with the usual norm, and let
$\ell_1$ denote the usual dual space of absolutely
summable sequences. Recall $\ell_1 = (c_0)^*$ and
$\ell_\infty = (\ell_1)^*$. Given Banach spaces $X,Y$ we
denote by $B(X,Y)$ the space of all bounded operators
$u\colon \ X\to Y$ equipped with the usual operator norm.
We will always identify an operator on a sequence space
with a matrix in the usual way. We will denote by $A_0$
(resp. $A_1$) the space of all complex matrices $(a_{ij})$
such that $$\sup_i \sum_j |a_{ij}| < \infty \quad
\left(\hbox{resp.} \sup_j \sum_i |a_{ij}|<\infty\right)$$
equipped with the  norm $$\eqalign{\|(a_{ij})\|_{A_0} &=
\sup_i \sum_j |a_{ij}|.\cr \bigg(\hbox{resp.} \quad
\|(a_{ij})\|_{A_1} &= \sup_j \sum_i |a_{ij}|.\bigg)}$$ We
will need to work with complex spaces, so we recall that
if $E$ is a real Banach lattice its complexification
$E+iE$ can be naturally equipped with a norm so that for
all $f = a+ib$ in $E+iE$ we have $\|a+ib\| = \|(|a|^2 +
|b|^2)^{1/2}\|_E$. We will call the resulting complex
Banach space a complex Banach lattice.

Let $E,F$ be real or complex Banach lattices. We will
denote by $B_r(E,F)$ the space of all operators $u\colon \
E\to F$ for which there is a constant $C$ such that for
all finite sequences $(x_i)_{i\le m}$ in $E$ we have
$$\|\sup_{i\le m}|u(x_i)|\, \|_F \le C\|\sup_{i\le m}|x_i|\,
\|_E.\leqno (1)$$ We will denote by $\|u\|_r$ the smallest
constant $C$ for which this holds, i.e. we set $\|u\|_r =
\inf\{C\}$. It is known, under the assumption that
 $F$ is Dedekind
complete in the sense of [MN] (sometimes also called order
complete), that  in the real (resp. complex) case, every
regular operator $u\colon \ E\to F$ is of the form $u =
u_+-u_-$ (resp. $u = a_+ - a_- + i(b_+ - b_-)$) where
$u_+, u_-$ (resp. $a_+, a_-, b_+, b_-$) are bounded
positive operators from $E$ to $F$, cf. [MN] or [S] p.233.
Since the converse is obvious this gives a very clear
description of the space $B_r(E,F)$. Here of course
``positive'' means positivity preserving. Let $E,F$ be
real (resp. complex) Banach lattices. Under the same
assumption on $F$ (cf. e.g. [MN] p.27) it is known that
$B_r(E,F)$ equipped with the usual ordering is a real
(resp. complex) Banach lattice in such a way that we have
$$\forall \ T\in B_r(E,F)\qquad \|T\|_r = \|\, |T|\,
\|_{B(E,F)}.$$ We refer to [MN] for more information on
the spaces $B_r(E,F)$. We will only use the following
elementary particular cases.

\n If $A$ is in $B_r(\ell_p,\ell_p)$ with associated matrix
$(a_{ij})$, let us denote by $|A|$ the operator admitting
$(|a_{ij}|)$ as its associated matrix. Then we have
$$\|A\|_r = \|\, |A|\, \|_{B(\ell_p,\ell_p)}. \leqno (2)$$
By a well known result (going back, I believe, to Grothendieck) for any
measure spaces $(\Omega,\mu)$, $(\Omega',\mu')$ we have an isometric
identity
$$B(L_1(\mu), L_1(\mu')) = B_r(L_1(\mu), L_1(\mu')).\leqno (3)$$
On the other hand, we have trivially (isometrically)
$$B(L_\infty(\mu), L_\infty(\mu')) = B_r(L_\infty(\mu),
L_\infty(\mu')) .\leqno (4)$$

The next result is known to many people
 in some slightly different form (in particular see [W]),
I believe that our formulation is useful and hope to
demonstrate this in the rest of this note. Motivated by
the results in [HeP], I suspected that this result was
known and I asked F.~Lust-Piquard whether she knew
 a reference for this, she did not but she immediately
showed me the following proof.

\proclaim Theorem 1. Let $A_0,A_1$ be as above. For any fixed integer $n$,
let $A^n_0 \subset A_0$, $A^n_1 \subset A_1$ be the subspace of all
matrices
$(a_{ij})$ which are supported by the upper left $n\times
n$ corner, so that the elements of $A^n_0$ or $A^n_1$ can
be viewed as $n\times n$ matrices.
Note the elementary identifications
$$\eqalign{A^n_0 &= B(\ell^n_\infty, \ell^n_\infty) \quad \hbox{and}\quad
A^n_1 = B(\ell^n_1, \ell^n_1)\cr
A_0 &= B(c_0, \ell_\infty) \quad \hbox{and}\quad A_1 = B(\ell_1,\ell_1).}$$
We have then for all $0 < \theta< 1$ the following isometric identities
where $p = 1/\theta$:\medskip
\item{(i)} $(A^n_0, A^n_1)_\theta = B_r(\ell^n_p, \ell^n_p)$.
\item{(ii)} $A_\theta = (A_0,A_1)^\theta = B_r(\ell_p,\ell_p)$.\medskip

\n {\bf Remark.} Note that $A_0$ and $A_1$ (resp. $A^n_0$ and $A^n_1$) are
isometric as Banach spaces. This is a special case of the fact that the
transposition induces an isometric isomorphism between the spaces
$B(X,Y^*)$ and $B(Y,X^*)$ when $X$ and $Y$ are Banach spaces.\medskip

\n {\bf Proof of Theorem 1.} (The main point was shown to me by
F.~Lust-Piquard.) We will prove (i) only. The second part (ii) follows
easily from (i) by a weak-$*$ compactness argument which we leave to the
reader.

\n Now let $(X_0,X_1)$ and $(Y_0,Y_1)$ be compatible
couples of finite dimensional complex Banach spaces and let
$X_\theta = (X_0,X_1)_\theta$, $Y_\theta =
(Y_0,Y_1)_\theta$.	It is well known and easy to check from
the definitions that we have a norm 1 inclusion
$$(B(X_0,Y_0), B(X_1,Y_1))_\theta \subset B(X_\theta,
Y_\theta).$$ Applying  this to the spaces $X_0 = Y_0 =
\ell^n_\infty(\ell^m_\infty)$ and $X_1 = Y_1 =
\ell^n_1(\ell^m_\infty)$ with $m$ an arbitrary integer, we
obtain the norm 1 inclusion $$(A^n_0, A^n_1)_\theta \to
B_r(\ell^n_p,\ell^n_p).$$ To show the converse, consider
the spaces $$B^n_0 = A^{n*}_0 \quad \hbox{and} \quad B^n_1
= A^{n*}_1.$$ For any $n\times n$ matrix $(b_{ij})$ we
have $$\|(b_{ij})\|_{B^n_0} = \sum^n_{i=1} \sup_{j\le n}
|b_{ij}| \quad \hbox{and} \quad \|(b_{ij})\| _{B^n_1} =
\sum^n_{j=1} \sup_{i\le n} |b_{ij}|.$$ We claim that for
all $(b_{ij})$ in the unit ball of $B^n_\theta
=(B^n_0)^{1-\theta} (B^n_1)^\theta = (B^n_0,B^n_1)_\theta$
we have $$\forall \ A \in B_r(\ell^n_p,\ell^n_p)\qquad
\sum_{i,j} |a_{ij}b_{ij}| \le \|A\|_r. \leqno (5)$$
With this claim, we conclude easily since (5) yields a
norm one inclusion $B_r(\ell^n_p, \ell^n_p)\subset
(B^n_\theta)^*$, and $(B^n_\theta)^* = (A^{n*}_0,
A^{n*}_1)^*_\theta = A^n_\theta$. Therefore it suffices to
prove the claim (5). Since $(b_{ij})$ is assumed in the
unit ball of $(B^n_0)^{1-\theta} (B^n_1)^\theta$ there are
$n\times n$ matrices $(b^0_{ij})$ and $(b^1_{ij})$ such
that $|b_{ij}|  = |b^0_{ij}| \cdot |b^1_{ij}|$ and such
that $\sum\limits_i \sup\limits_j |b^0_{ij}|^{p'} \le 1$,
$\sum\limits_j \sup\limits_i |b^1_{ij}|^p \le 1$, with $p
= 1/\theta$ and $p'  = 1/(1-\theta)$.

Now let $\beta_i = \sup\limits_j|b^0_{ij}|$ and $\alpha_j  = \sup\limits_i
|b^1_{ij}|$. Then
$$\eqalignno{\sum |a_{ij}b_{ij}| &\le \sum_{i,j} \beta_i|a_{ij}|\alpha_j\cr
\noalign{\hbox{hence by (2)}}
&\le \|A\|_r.}$$
This proves our claim and concludes the proof.\qed

It is then routine to deduce the following extension.

\proclaim Corollary 2. Let $(\Omega,\mu)$ and $(\Omega',\mu')$ be arbitrary
measure spaces. Consider the couple
$$X_0 = B(L_\infty(\mu), L_\infty(\mu')), \quad X_1 = B(L_1(\mu),
L_1(\mu')).$$
We will identify (for the purpose of interpolation) elements in $X_0$ or
$X_1$ with linear operators from the space of integrable
step functions into $L_1(\mu') + L_\infty(\mu')$. We have
then isometrically
$$X^{1-\theta}_0 X_1^{\theta} = (X_0,X_1)^\theta = B_r(L_p(\mu),
L_p(\mu')).\leqno (6)$$

\n {\bf Remark.} Recalling (3) and (4), we can rewrite (6) as follows
$$(B_r(L_\infty(\mu), L_\infty(\mu')), B_r(L_1(\mu), L_1(\mu')))^\theta =
B_r(L_p(\mu), L_p(\mu')).$$
By [Be], it follows that the space $(X_0,X_1)_\theta$ coincides with the
closure in $B_r(L_p(\mu), L_p(\mu'))$ of the subspace of all the operators
which are simultaneously bounded from $L_1(\mu)$ to $L_1(\mu')$ and from
$L_\infty(\mu)$ to $L_\infty(\mu')$.

We will now consider operators defined on a subspace
$S$ of a Banach lattice $E$ and taking values in a Banach
lattice $F$. Let $u:S\to F$ be such an operator. We will
again say that
$u$ is {\it regular} if there is a constant $C$ such that
$u$ satisfies (1) for all finite sequences
$x_1,...,x_m$ in $E$. We again denote by $\|u\|_r $ the
smallest constant $C$ for which this holds. Clearly the
restriction to $S$ of a regular operator defined on $E$ is
regular. Conversely, in general a regular operator
on $S$ is not necessarily the restriction of a regular
operator on $E$: for instance if $E$ is $L_1$, if $S$
is the closed span of a sequence of standard independent
Gaussian random variables and if $u:S\to L_2$ is the
natural inclusion map, then $u$ is regular (this is a
well known result of Fernique, see e.g. [LeT] p. 60) but
does not extend to any bounded map from $L_1$ into $L_2$
since by a weakening of Grothendieck's theorem (cf. [P4]
p. 57), the identity of $S$  would then be $2$-absolutely
summing, which is absurd since $S$ is infinite dimensional
(cf. e.g. [P4] p. 14).

Nevertheless, it
turns out that in several interesting cases, conversely
every regular operator on $S$ is the restriction of a
regular operator on $E$ with the same regular norm.  In
particular, the next statement is  an extension theorem for
regular operators which generalizes a result due to
M.~L\'evy [L\'e] in the case $p=1$. We will prove

\proclaim Theorem 3. Let $1 \le p \le \infty$. Let $(\Omega,\mu)$,
$(\Omega',\mu')$ be arbitrary measure spaces. Let $S\subset L_p(\mu)$ be
any closed subspace. Then every regular operator   
$u\colon \ S\to L_p(\mu')$ admits a regular extension
$\tilde u\colon \ L_p(\mu) \to L_p(\mu')$ such that
$\|\tilde u\|_r = \|u\|_r$.

Actually this will be a consequence of the following more general result.
(We refer the reader to [LT] for the notions of
$p$-convexity and $p$-concavity.)

  \proclaim Theorem 4. Let
$L,\Lambda$ be Banach lattices and let $S\subset \Lambda$
be a closed subspace. Assume that $L$ is a dual space, or
merely that there is a regular projection $P\colon \
L^{**}\to L$ with $\|P\|_r \le 1$. Assume moreover that
for some $1\le p \le \infty$ $\Lambda$ is $p$-convex and
$L$ $p$-concave. Then every regular operator $u\colon \
S\to L$ extends to a regular operator $\tilde u\colon \
\Lambda\to L$ with $\|\tilde u\|_r = \|u\|_r$.

\n {\bf Remark.} Note that by known results (cf. [K] or
[LT]) in the above situation every positive operator
$u\colon \ \Lambda\to L$ factors through an $L_p$-space,
i.e. there is a measure space $(\Omega,\mu)$ and operators
$B\colon\ \Lambda \to L_p(\mu)$ and $A\colon \ L_p(\mu)\to
L$ such that $U = AB$ and $\|A\|\cdot \|B\| = \|u\|$.
Actually for this conclusion to hold, it suffices to
assume that $u$ can be written as the composition of
first a $p-$convex operator with constant $\le 1$ followed
by a $p$-concave operator with constant $\le 1$. Therefore,
since every regular operator with regular norm $\le 1$ on a
$p$-convex Banach lattice clearly is itself $p$-convex
with constant $\le 1$, 
every {\it regular\/} $u\colon \ \Lambda\to L$ 
factors through an $L_p$-space with factorization
constant at most $1$. Actually it is easy to modify
Krivine's argument to prove that, in the same situation
as in Theorem 4, 
every {\it regular\/} $u\colon \ \Lambda\to L$ 
can be written as $u = AB$ as above but with $A,B$
regular and such that $\|A\|_r \|B\|_r =\|u\|_r$. \medskip

\n {\bf Proof of Theorem 4.} By a standard ultraproduct argument it is
enough to consider the case when $L$ is finite dimensional with an
unconditional basis $(e_1,\ldots, e_n)$. As usual in extension problems, we
will use the Hahn-Banach theorem. We need to introduce a Banach space $X$
such that $X^* = B_r(\Lambda, L)$. The space $X$ is defined as the tensor
product $L^*\otimes \Lambda$ equipped with the following norm, for all $v =
\sum^n_1 \alpha_ke^*_k \otimes x_k$ with $\alpha_i$ scalar and $x_i \in
\Lambda$ we define
$$\|v\|_X = \inf\left\{ \left\|\sum\nolimits^n_1 \alpha_ie_i\right\|_{L^*}
\|\sup_{i\le n} |x_i|\,\|_\Lambda\right\}.$$
The only assumption needed for our extension theorem is that $\|\quad\|_X$
is a norm (see the remark below). This follows from the $p'$-convexity of
$L^*$ and the $p$-convexity of $\Lambda$.

To check this we assume as we may that $\Lambda$ is
included in a space of measurable functions $L_0(\mu)$ on
some measure space. Let $Y_0$ be the space of $n$-tuples
of measurable functions $y_1,\ldots, y_n$ in
$L_\infty(\mu)$ equipped with the norm $$\|(y_i)\|_{Y_0} =
\left\|\sum^n_1 \|y_i\|^{1\over p'}_{_\infty}
e^*_i\right\|^{p'}_{L^*}.$$ That this is indeed a norm
follows from the $p'$-convexity of $L^*$.

Let $Y_1$ be the space of $n$-tuples of measurable functions
$y_1,\ldots,y_n$ is $L_0(\mu)$ such that $|y_i|^{1\over p} \in \Lambda$
equipped with the norm
$$\|(y_i)\|_{Y_1} = \|\sup_{i\le n} |y_i|^{1\over p}\|^p_\Lambda.$$
Again this is a norm by the $p$-convexity of $\Lambda$. But now if we
consider the unit ball of the space
$$Y^{1-\theta}_0 Y^\theta_1 \quad \hbox{with}\quad \theta = {1\over p}$$
we find exactly the set $C$. This shows that $C$ is convex as claimed
above. We will now check that $X^* = B_r(\Lambda, L)$ isometrically.

Consider $u\colon \ \Lambda\to L$. We have
$$\|u\|_r = \sup\left\{\left\| \sum^n_{k=1} \sup_{i\le m} |\langle u(x_i),
e^*_k\rangle| e_k\right\|_L \right\}\leqno (7)$$
where the supremum runs over all $m$ and all $m$-tuples $(x_1,\ldots, x_m)$
in $\Lambda$ such that\break $\|\sup\limits_{i\le m} |x_i|\, \|_\Lambda \le
1$.
Let us denote by $\beta$ the unit ball of $L^*$. Then (7) can be rewritten
$$\|u\|_r = \sup\left\{\left| \sum^n_{k=1} \alpha_k \langle u(x_{i_k}),
e^*_k\rangle \right|\right\}\leqno (8)$$
where the supremum runs over all integers $m$, all choices $i_1,\ldots,
i_n$ in $\{1,\ldots, m\}$, all elements $\alpha =  \sum\limits^n_1
\alpha_ke^*_k$ in $\beta$ and all $m$-tuples $x_1,\ldots, x_m$ in $\Lambda$
with $\|\sup|x_i|\, \|_\Lambda \le 1$. But for such elements clearly $v =
\sum\limits^n_{k=1} \alpha_k e^*_k\otimes x_{i_k}$ is in the set $C$ which
is the unit ball of $X$, hence (8) yields
$$\|u\|_r = \sup\{|\langle u,v\rangle|\mid v\in B_X\}.$$
This proves the announced claim that $X^* = B_r(\Lambda, L)$ isometrically.

We can then complete the proof by a well known application of the
Hahn-Banach theorem.

Consider the subspace $M\subset X$ formed by all the $v = \sum^n_1 \alpha_k
e^*_k \otimes x_k$ such that $\alpha_kx_k \in S$ for all $k=1,\ldots, n$.
If $u\colon\ S \to L$ is regular we clearly have for all $v$ in $M$
$$|\langle u,v\rangle| = \left|\sum^n_1 \langle u(\alpha_kx_k),
e^*_k\rangle\right| \le \|u\|_r \|v\|_X$$
hence we can find a Hahn-Banach extension of the linear form $v\in M\to
\langle u,v\rangle$ defined on the whole of $X$ and still with norm $\le
\|u\|_r$. Clearly we can write the extension in the form $v\in X\to \langle
\tilde u, v\rangle$ for some operator $\tilde u\colon \ \Lambda \to L$ and
since $\sup\limits_{\|v\|_X\le 1} |\langle \tilde u,v\rangle| \le\|u\|_r$,
we have $\|\tilde u\|_r \le \|u\|_r$ as announced.\qed \medskip

\n {\bf Remark.} Assume again $L$ finite dimensional as above. The
assumption ``$\Lambda$ $p$-convex, $L$ $p'$-concave'' can be replaced by
the property that in $L^*\otimes \Lambda$ the set
$$\eqalign{C = &\left\{\sum\nolimits^n_1 \alpha_ie^*_i \otimes x_i\mid
\alpha_i \in {\bf C} \quad x_i \in \Lambda\right.\cr
&\left.\left\|\sum\nolimits^n_1 \alpha_ie^*_i\right\|_{L^*} \le 1, \quad
\|\sup|x_i|\|_\Lambda \le 1\right\}}$$
is a convex set.

As the preceding proof shows this is true is $L^*$ is $p'$-convex and
$\Lambda$ $p$-convex. However, it is clearly true also in other cases. For
instance if $L^*  = \ell^n_\infty$ then $C$ is just the unit ball of
$\Lambda(\ell^n_\infty)$ which is clearly convex for all $\Lambda$. Moreover
if $\Lambda = L_\infty(\mu)$ for some measure $\mu$ then $C$ is the unit
ball of $L^*(L_\infty(\mu))$ which is convex for all $L$. More generally,
what we really use (and which is then equivalent
 to the extension theorem, by a reasoning well known to
many Banach space specialists) is that the closed convex
hull of the set $C$, satisfies $$\overline{\rm conv}(C)
\cap M = \overline{\rm conv}(C\cap M),$$ where $M$ denotes
as above the subspace $M = L^* \otimes S \subset L^*
\otimes \Lambda$.

We now give some applications to $H^p$-spaces, mainly
motivated by our paper [P2].
 Let $1\le p \le \infty$. Let $H^p$ be the
usual $H^p$-space of functions on the torus $\bf T$ equipped with its
normalized Haar measure $m(dt) = {dt\over {2\pi}}$. We denote simply $L_p =
L_p({\bf T}, m)$. Given a finite dimensional normed space $E$ we denote
$$H^p(E) = \{f\in L_p(m;E)\mid \hat f(n) = 0\quad
 \forall\ n< 0\}.$$
By a result of P. Jones [J] (see also [BX, P1, X] for a
discussion of the vector valued case) we have
isomorphically and with isomorphism constants independent
of $E$ $$H^p(E) = (H^\infty(E), H^1(E))_\theta \quad
\hbox{if}\quad \theta =1/p. \leqno \hbox{(9)}$$ More
precisely, there is a constant $C_p$ such that for all $f$
in $H^p(E)$ we have  $$\|f\|_{(H^\infty(E), H^1(E))_\theta}
\le C_p \|f\|_{H^p(E)}. \leqno(10)$$
 We will prove the following
extension of Corollary~2.

\proclaim Theorem 5. Let $(\Omega,\mu)$ be an arbitrary measure space. Let
$$B_0 = B_r(H^\infty, L_\infty(\mu)) \quad B_1 = B_r(H^1, L_1(\mu)).$$
Then (isomorphically) $(B_0,B_1)^\theta = B_r(H^p,
L_p(\mu))$ with $\theta = 1/p$.

\pf By Theorem 4, if $u\colon \ H^p \to L_p(\mu)$ is such that $\|u\|_r
<1$, then $\exists\ \tilde u\colon \ L_p\to L_p(\mu)$ extending $u$ such
that $\|\tilde u\|_r < 1$.
 By Corollary~2, $\tilde u$ is of norm $<1$ in the space
$(B_r (L_\infty, L_\infty(\mu))$, $B_r(L_1, L_1(\mu)))^\theta$ hence by
restriction $u$ is of norm $<1$ in $(B_0,B_1)^\theta$. Conversely, assume
that $u$ is in the unit ball of $(B_0,B_1)^\theta$. Consider then $f$ in
the unit ball of $H^p(\ell^n_\infty)$, or equivalently consider an
$n$-tuple $(f_1,\ldots, f_1)$ in $H^p$ such that $\int \sup\limits_{k\le n}
|f_k|^p dm \le 1$. By P.~Jones's theorem (10) we have
$$\|f\|_{(H^\infty(\ell^n_\infty), H^1(\ell^n_\infty))_\theta} \le
C_p,$$ hence since $\|u\|_{(B_0,B_1)^\theta} \le 1$ by
assumption, it is easy to deduce
$$\|u(f)\|_{(L_\infty(\ell^n_\infty), L_1(\ell^n_\infty))^\theta} \le C_p$$
or equivalently since $L_p(\ell^n_\infty) = (L_\infty(\ell^n_\infty),
L_1(\ell^n_\infty))^\theta$
$$\int \sup_{k\le n} |u(f_k)|^p dm \le C_p^p.$$
By homogeneity we conclude that
$\|u\|_{B_r(H^p, L_p(\mu))} \le C_p$.\qed

\n {\bf Remark.} Once again by [Be], the space $(B_0,B_1)_\theta$ coincides
with the closure in $B_r(H^p, L_p(\mu))$ of the operators which are
simultaneously regular from $H^1$ to $L_1(\mu)$ and from $H^\infty$ to
$L_\infty(\mu)$.

\n {\bf Remarks.}  \item {(i)} For a version of Theorem 1
and Corollary 2 in the case of noncommutative $L_p$-spaces,
we refer the reader to [P3].

 \item{(ii)} Of course Theorem 5 and its proof remain valid
with the couple  $(H^\infty,H^1)$ replaced by any couple of
subspaces of $(L^\infty,L^1)$ for which (10) holds.

\vskip12pt
\centerline{\bf References}
\vskip12pt
\item{[B]} J. Bergh. On the relation between the two
complex methods of interpolation. Indiana Univ. Math.
Journal 28 (1979) 775-777.
 
 \item{[BL]} J. Bergh and J. L\"ofstr\"om. Interpolation spaces. An
 introduction. Springer Verlag, New York. 1976.
 
\item{[BX]} O. Blasco and Q. Xu. Interpolation between
vector valued Hardy spaces, J. Funct. Anal. 102 (1991)
331-359.

 \item{[Ca]} A. Calder\'on. Intermediate spaces and interpolation, the
 complex method. Studia Math. 24 (1964) 113-190.
 
\item{[HP]} U. Haagerup and G. Pisier. Factorization of
analytic functions
 with values in non-commutative $L^1$-spaces and applications. Canadian J.
 Math. 41 (1989) 882-906.

\item{[HeP]} A. Hess and G. Pisier, On the
$K_t$-functional for the couple $B(L_1,L_1), B(L_\infty,
L_\infty)$).

\item{[K]} J.L. Krivine. Th\'eor\`emes de factorisation dans
les espaces de Banach r\'eticul\'es. S\'eminaire
Maurey-Schwartz 73/74, Expos\'e 22, Ecole Polytechnique,
Paris.

\item {[L\'e]} M. L\'evy. Prolongement d'un op\'erateur
d'un sous-espace de $L^1(\mu)$ dans $L^1(\nu)$. S\'eminaire
d'Analyse Fonctionnelle 1979-1980. Expos\'e 5. Ecole
Polytechnique.Palaiseau.

\item{[LeT]} M. Ledoux and M. Talagrand. Probability in
Banach spaces. Springer-Verlag 1991.

\item{[LT]} J. Lindendrauss and L. Tzafriri. Classical
Banach spaces II, Function spaces, Springer-Verlag, 1979.

\item{[MN]} P. Meyer-Nieberg. Banach Lattices.
Universitext, Springer-Verlag,  1991.

\item{[P1]} G. Pisier. Interpolation of $H^p$-spaces
and noncommutative generalizations I. Pacific J. Math. 155
(1992) 341-368.

\item{[P2]}$\underline{\hskip1.5in}$. Interpolation of $H^p$-spaces
and noncommutative generalizations II. Revista Mat.
Iberoamericana. To appear. 

\item{[P3]}$\underline{\hskip1.5in}$. The Operator
Hilbert space $OH$, Complex Interpolation and Tensor
Norms. To appear.

\item {[P4]} $\underline{\hskip1.5in}$. Factorization of linear
operators and the Geometry of Banach spaces.  CBMS
(Regional conferences of the A.M.S.)    60, (1986),
Reprinted with corrections 1987.

\item{[S]} H.H. Schaefer. Banach lattices and positive
operators. Springer-Verlag, Berlin Heidelberg New-York,
1974.

\item{[W]} L. Weiss. Integral operators and changes of
density. Indiana University Math. Journal 31 (1982) 83-96.

\item{[X]} Q. Xu. Notes on interpolation of Hardy spaces.
Ann. Inst. Fourier. To appear.

\vskip12pt

\end